%% file: outlier_isit_final.tex
\documentclass[conference]{IEEEtran}
\IEEEoverridecommandlockouts
\input{preamble}

\usepackage{subfigure,epstopdf,bbm} 
\usepackage[utf8]{inputenc} 
\usepackage[T1]{fontenc}    
\usepackage{url,bbm}            
\usepackage{booktabs}       
\usepackage{amsfonts}       
\usepackage{nicefrac}       
\usepackage{microtype}      
\usepackage{cite}
\usepackage{subfigure,transparent,color,graphicx} 

\allowdisplaybreaks[2]
\def\BibTeX{{\rm B\kern-.05em{\sc i\kern-.025em b}\kern-.08em T\kern-.1667em\lower.7ex\hbox{E}\kern-.125emX}}

\begin{document}
\title{Asymptotics for Outlier Hypothesis Testing\\
{\thanks{This work was partially supported by a National Key Research and Development Program of China under Grant 2020YFB1804800, and by a grant from the US Army Research Office (ARO W911NF-15-1-0479).}}
}

\author{  
 \IEEEauthorblockN{Lin Zhou}
  \IEEEauthorblockA{School of Cyber Science and Technology\\
Beihang University\\
    Email: \url{lzhou@buaa.edu.cn}}
       \and
  \IEEEauthorblockN{Yun Wei}
  \IEEEauthorblockA{Department of Statistical Science\\
       Duke University \\
   Email: \url{yun.wei@duke.edu}}
  \and
  \IEEEauthorblockN{Alfred Hero III}
  \IEEEauthorblockA{Department of EECS\\
       University of Michigan \\
   Email: \url{hero@eecs.umich.edu}
   }
}
\maketitle

\begin{abstract}
 We revisit the outlier hypothesis testing framework of Li \emph{et al.} (TIT 2014) and derive fundamental limits for the optimal test under the generalized Neyman-Pearson criterion. In outlier hypothesis testing, one is given multiple observed sequences, where most sequences are generated i.i.d. from a nominal distribution. The task is to discern the set of outlying sequences that are generated according to anomalous distributions. The nominal and anomalous distributions are \emph{unknown}. We consider the case of multiple outlying sequences where the number of outlying sequences is unknown and each outlying sequence can follow a different anomalous distribution. Under this setting, we study the tradeoff among the probabilities of misclassification error, false alarm and false reject. Specifically, we propose a threshold-based test that ensures exponential decay of misclassification error and false alarm probabilities. We study two constraints on the false reject probability, with one constraint being that it is a non-vanishing constant and the other being that it has an exponential decay rate. For both cases, we derive bounds on the false reject probability, as a function of the threshold, for each tuple of nominal and anomalous distributions.
\end{abstract}


\section{Introduction}

Motivated by practical applications in anomaly detection~\cite{chandola2009anomaly}, we revisit the outlier hypothesis testing (OHT) problem studied in~\cite{li2014}. In the OHT problem, one is given $M$ sequences and asked to discern the set of outlying sequences which are generated from an \emph{unknown} anomalous distribution that is different from an \emph{unknown} nominal distribution from which the rest majority of the sequences are generated from. We consider the case where the outlying sequence might \emph{not} be present and derive the performance tradeoff between the probabilities of misclassification error, false alarm and false reject for a threshold-based test. Furthermore, we show that our test is optimal under the generalized Neyman-Pearson criterion~\cite{gutman1989asymptotically} for both a second-order asymptotic regime and a large deviations regime. Our second-order asymptotic result 
approximates the finite sample performance of our test. Throughout the paper, we assume that the sequences have a finite alphabet.

We assume that the number of outlying sequence (outlier) is \emph{unknown} and each outlying sequence can be drawn from a different anomalous distribution. When the number of outliers is known, Li \emph{et al.}~\cite[Theorem 10]{li2014} derived an achievability decay rate of the error probabilities under each hypothesis and showed asymptotic optimality of their result when the number of the sequences $M$ tends to infinity, when the lengths of sequences $n$ tend to infinity and when all the outlying sequences are generated from the same anomalous distribution. Furthermore, when the number of outliers is unknown and when each outlier is generated from the same anomalous distribution, Li \emph{et al.}~\cite[Theorem 10]{li2014} showed that when the null hypothesis is not taken into account, a generalized likelihood ratio test is exponentially consistent. However, the authors of \cite{li2014} did not characterize the exponent explicitly. One might wonder whether it is possible to characterize the performance of a test when the number of outliers is unknown and when each outlier can be generated from a different anomalous distribution. We answer this question affirmatively by proposing a threshold-based test, characterizing its performance explicitly and proving its asymptotic optimality under the generalized Neyman-Pearson criterion~\cite{gutman1989asymptotically}.

\subsection{Main Contributions}
\vspace{-0.2em}

Our main contributions are two fold. Firstly, we propose a threshold-based test in \eqref{test:Toutlier} that is ignorant of the nominal and anomalous distributions and the number of outliers. Secondly, we analyze the tradeoff among probabilities of misclassification error, false reject and false alarm for our test. Specifically, under each tuple of unknown nominal and anomalous distributions, we show that our test ensures that both false alarm and misclassification error probabilities decay exponentially fast and we bound the false reject probability as a function of the threshold in two regimes. In the first regime named second-order asymptotics, we derive bounds on the false reject probability when the number of samples is finite and show that asymptotically when the lengths of the observed sequences tend to infinity, the false reject probability is upper bounded by a constant $\varepsilon\in(0,1]$. Furthermore, we also study the second regime named large deviations where asymptotically we derive the exponential decay rate of the false reject probability as a function of the threshold in our test. We establish that, as long as the nominal and anomalous distributions are far in a given distance measure that generalizes Jensen-Shannon divergence~\cite{lin1991divergence},  our test is exponentially consistent by ensuring that all three error probabilities decay exponentially. For both regimes, we show that our test is optimal under the generalized Neyman-Pearson criterion.

\subsection{Related Works}
The most closely related work to ours is that of \cite{li2014,zhou2021itw}. In \cite{li2014}, the authors formulated the outlier hypothesis testing problem, and derived optimal results under constraints on the number of observed sequences, the length of observed sequences and the number of anomalous distributions. In our previous work~\cite{zhou2021itw}, we revisited the case of at most one outlier in \cite{li2014} and derived bounds on the false reject probability under any pair of nominal and anomalous distributions for the generalized Neyman-Pearson setting where misclassification error and false alarm probabilities are constrained to decay exponentially fast for all pairs of  distributions. However, such a setting twists the achievability and converse proofs with the test design and is thus confusing\footnote{We thank the reviewers of our IT version~\cite{zhou2020second} for these comments.}. In this paper, we generalize \cite{zhou2021itw} to the case of multiple outliers and solve the above problem by presenting the test design, the achievability and converse results separately in different subsections. Our results, when specialized to the case of at most one outlier, present those in \cite{zhou2021itw} in a clearer manner. Other non-exhausted related work on outlier hypothesis testing includes \cite{bu2019linear,li2017universal,cohen2015active,Zou2017TSP}.

Since our proof technique is inspired by asymptotic statistical classification theory, we also mention a few works in this domain. In \cite{gutman1989asymptotically}, the author studied a binary sequence classification problem and showed that a certain test using empirical distributions is asymptotically optimal with exponentially decreasing misclassification probabilities. The result in \cite{gutman1989asymptotically} was generalized to classification of multiple sequences in~\cite{unnikrishnan2015asymptotically} and to distributed detection in~\cite{he2019distributed}. Finally, a finite sample analysis for the setting of~\cite{gutman1989asymptotically} was provided in~\cite{zhou2018binary}.

\section{Problem Formulation and Test Design}
\label{sec:Tout}

\subsection{Problem Formulation}
Given any $M\in\bbN$, let $T:=\lceil \frac{M}{2}-1\rceil$. For any integer $t\in[T]$, let $\calS_t$ denote the set of all subsets of $[M]$ whose cardinality (size) is $t$, i.e.,
\begin{align}
\calS_t:=\{\calB\subseteq[M]:~|\calB|=t\}\label{def:calSt}.
\end{align}
Then, define the union of subsets $\calS$ over $t\in[T]$ as follows:
\begin{align}
\calS:=\bigcup_{t\in[T]}\calS_t\label{def:calS}.
\end{align}

In the outlier hypothesis testing problem with at most $T$ outliers, the task is to decide whether there are outliers among $M$ observed sequences $\bX^n=(X_1^n,\ldots,X_M^n)$ and identify the set of outlying sequences if any exist. We assume that each outlying sequence is generated i.i.d. according to a possibly different anomalous distribution. Specifically, let $\bP_T:=(P_{\rmA,1},\ldots,P_{\rmA,T})$ be a collection of $T$ anomalous distributions that are different from the nominal distribution $P_\rmN$, all defined on the finite alphabet $\calX$ with the same support. Furthermore, for any $\calB\in\calS$, let $\bP_\calB$ denote the collection of distributions $(P_1,\ldots,P_{|\calB|})$. When $\calB\in\calS$ denotes the index of the outlying sequences, for any $l\in\calB$, $X_l^n$ is generated i.i.d. from $P_{\rmA,\jmath_\calB(t)}$, where $\jmath_\calB$ denotes an ordered mapping from $\calB$ to $[|\calB|]$ such that for each $i\in\calB$, $\jmath_\calB(i):=j$ if $i$ is the $j$-th smallest element in $\calB$. For example, when $M=10$, $\calB=\{2,3,6\}$ and $\bP_\calB=(P_{\rmA,1},P_{\rmA,2},P_{\rmA,3})$, then the second sequence $X_2^n$ is generated i.i.d. from $P_{\rmA,1}$, the third sequence $X_3^n$ is generated i.i.d. from $P_{\rmA,2}$ and the $6$-th sequence $X_6^n$ is generated i.i.d. from $P_{\rmA,3}$ while all other sequences are generated i.i.d. from the unknown nominal distribution $P_\rmN$.

Since the exact number of outliers is \emph{unknown}, there are in total $|\calS|+1=\sum_{t\in[T]}{M \choose t}+1$ possible configurations of outlying sequences. Formally, the task is to design a test $\phi_n:\calX^{Mn}\to \{\{\rmH_\calB\}_{\calB\in\calS},\rmH_\rmr\}$, ignorant of the distributions $(P_\rmT,\bP_T)$, to classify among the following hypotheses:
\begin{itemize}
\item $\rmH_\calB$ where $\calB\in\calS$: the set of outlying sequences are sequences $X_j^n$ with $j\in\calB$;
\item $\rmH_\rmr$: there is no outlying sequence.
\end{itemize}
The null hypothesis $\rmH_\rmr$ is introduced to model the case when there is no outlier among all $M$ observed sequences.

Given any test $\phi_n$, under any tuple of nominal and anomalous distributions $(P_\rmN,\bP_T)=(P_\rmN,P_{\rmA,1},\ldots,P_{\rmA,T})$, the performance of  $\phi_n$ is evaluated by the following probabilities of misclassification error, false reject and false alarm:
\begin{align}
\beta_{\calB}(\phi_n|P_\rmN,\bP_T)
&:=\bbP_{\calB}\{\phi_n(\bX^n)\notin\{\rmH_\calB,\rmH_\rmr\}\}\label{def:errorsetS},\\
\zeta_{\calB}(\phi_n|P_\rmN,\bP_T)
&:=\bbP_{\calB}\{\phi_n(\bX^n)=\rmH_\rmr\}\label{def:rejectsetS},\\
\rmP_{\mathrm{fa}}(\phi_n|P_\rmN,\bP_T)
&:=\bbP_\rmr\{\phi_n(\bX^n)\neq\rmH_\rmr\}\label{def:falarmM},
\end{align}
where $\calB\in\calS$ denotes the set of indices of outliers, and we define $\bbP_{\calB}(\cdot):=\Pr\{\cdot|\rmH_\calB\}$ where for $i\in\calM_\calS$, $X_i^n$ is generated i.i.d. from the nominal distribution $P_\rmN$ and for $i\in\calB$, $X_i^n$ is generated i.i.d. from an nominal distribution $P_{\rmA,\jmath_\calB(i)}$, finally we define $\bbP_\rmr(\cdot):=\Pr\{\cdot|\rmH_\rmr\}$, where all sequences are generated i.i.d. from the nominal distribution $P_\rmN$.

\subsection{A Threshold-Based Test}
We use a threshold-based test that takes the empirical distribution of each observed sequence as the input and outputs a decision among all hypotheses. We need the following definition to present our test. Given a sequence of distributions $\bQ=(Q_1,\ldots,Q_M)\in\calP(\calX)^M$ and each $\calB\in\calS$, define the following linear combination of KL divergence terms
\begin{align}
\rmG_\calB(\bQ)
&:=\sum_{t\in\calM_{\calB}}D\left(Q_t\bigg\|\frac{\sum_{l\in\calM_\calB}Q_l}{M-|\calB|}\right)\label{def:gB},
\end{align}
where $\calM_\calB$ is defined as the set of elements that are in $[M]$ but not in $\calB$, i.e., $\calM_\calB:=\{i\in[M]:~i\notin\calB\}$. Note that $\rmG_\calB(\bQ)$ is a homogeneous measure and equals zero if and only if $Q_j=Q$ for all $j\in\calM_\calB$ where $Q\in\calP(\calX)$ is arbitrary

Given $M$ observed sequences $\bx^n=(x_1^n,\ldots,x_M^n)$ and any positive real number $\lambda$, our test operates as follows:
\begin{align}
\label{test:Toutlier}
\Psi_n(\bx^n)
&:=\left\{
\begin{array}{cc}
\rmH_\calB&\mathrm{if~}\rmS_\calB(\bx^n)<\min_{\calC\in\calS_\calB}\rmS_\calC(\bx^n)\\
&\mathrm{and~}\min_{\calC\in\calS_\calB}\rmS_\calC(\bx^n)>\lambda,\\
\rmH_\rmr&\mathrm{otherwise},
\end{array}
\right.
\end{align}
where $\calS_\calB=\calS\setminus\{\calB\}$ and $\rmS_\calC(\cdot)$ is the following scoring function:
\begin{align}
\rmS_\calC(\bx^n)
&:=\rmG_\calC(\hatT_{x_1^n},\ldots,\hatT_{x_M^n})\label{def:scalC},
\end{align}
which measures the sum of KL divergence between the empirical distribution of each sequence $x_i^n$ with $i\notin\calB$ relative to the average of the empirical distribution of all sequences $x_j^n$ where $j\notin\calB$.

\subsection{Discussions}
In this subsection, we discuss the performance of the test in \eqref{test:Toutlier}. For this purpose, we need the following definitions. Given any $\calB\in\calS$ and any tuple of distributions $\bP_\calB=(P_\rmN,P_{\rmA,1},\ldots,P_{\rmA,|\calB|})\in(\calP(\calX))^{|\calB|+1}$, for any two sets $(\calB,\calC)\in\calS^2$, define the following mixture distribution
\begin{align}
\nn P_{\rm{Mix}}^{(\calB,\calC,P_\rmN,\bP_\calB)}(x)
&:=\frac{1}{M-|\calC|}\Big(\sum_{i\in(\calB\cap\calM_{\calC})}P_{\rmA,\jmath_\calB(i)}(x)\\*
&\qquad+\sum_{i\in(\calM_{\calB}\cap\calM_{\calC})}P_\rmN(x)\Big)\label{def:Pmix},
\end{align}
and define the following sum of KL divergences
\begin{align}
\nn &\mathrm{GD}(\calB,\calC,P_\rmN,\bP_\calB)
=\sum_{i\in(\calB\cap\calM_{\calC})}D\left(P_{\rmA,\jmath_\calB(i)}\|P_{\rm{Mix}}^{(\calB,\calC,P_\rmN,\bP_\calB)}\right)\\*
&\qquad\qquad+\sum_{i\in(\calM_{\calB}\cap\calM_{\calC})}D\left(P_{\rmN}\|P_{\rm{Mix}}^{(\calB,\calC,P_\rmN,\bP_\calB)}\right)\label{def:GDBC}.
\end{align}

We first provide intuition into why the above test should be asymptotically consistent. Assume that $\calB\in\calS$ denotes the set of indices of the outliers and for each $l\in\calB$, the outlier $x_l^n$ is generated i.i.d. from $P_{\rmA,\jmath_{\calB}(l)}$. As $n$ tends to infinity, for any sequence $x_i^n$ where $i\notin\calB$, the empirical distribution $\hatT_{x_i^n}$ tends to the nominal distribution $P_\rmN$ and for each $x_i^n$ where $i\in\calB$, the empirical distribution $\hatT_{x_i^n}$ tends to $P_{\rmA,\jmath_{\calB}(i)}$. Therefore, given any $\calC\in\calS_\calB$, the scoring function $\rmS_\calC(\bx^n)$ converges to $\mathrm{GD}(\calB,\calC,P_\rmN,\bP_\calB)$ and $\rmS_\calB(\bx^n)$ converges to zero. Note that $\mathrm{GD}(\calB,\calC,P_\rmN,\bP_\calB)>0$ for any $\bP_\calB$ where $P_{\rmA,j}\neq P_\rmN$ for all $j\in[|\calB|]$. Thus, asymptotically, the set of outliers $\calB$ can be identified if $\lambda<\min_{\calC\in\calS_\calB}\mathrm{GD}(\calB,\calC,P_\rmN,\bP_\calB)$. On the other hand, when there is no outlier, for each $\calB\in\calS$, the scoring function $\rmS_{\calB}(\bx^n)$ tends to zero and thus, with any positive threshold $\lambda$, the null hypothesis is decided. Therefore, the test in \eqref{test:Toutlier} is consistent asymptotically for any set of outlier indices $\calB\in\calS$ and for any tuple of distributions $\bP_\calB$ where $P_{\rmA,j}\neq P_\rmN$ for all $j\in[|\calB|]$ such that the threshold satisfies $\lambda<\min_{\calC\in\calS_\calB}\mathrm{GD}(\calB,\calC,P_\rmN,\bP_\calB)$.

We then discuss how the test in \eqref{test:Toutlier} deals with unknown number of outliers when $T\geq 2$. Given $M$ observed sequences $\bx^n$, we calculate the scoring functions $\rmG_\calB(\hatT_{x_1^n},\ldots,\hatT_{x_2^n})$ for all possible sets $\calB\subseteq\calS$. Note that each $\calB\subseteq(\calS\setminus\emptyset)$ denotes a possible set of indices of outlying sequences and $\calB=\emptyset$ corresponds to the null hypothesis that no outlier appears. To determine the set of outliers, using the scoring function for all possible ${M\choose T}+1$ cases, we run the test in \eqref{test:Toutlier} that compares each scoring function with the threshold $\lambda$. In other words, the test \eqref{test:Toutlier} checks all possibilities of outliers to make a decision and its complexity increases exponentially with $T$.

Finally, we remark that the statistic in \eqref{def:gB} was also used in \cite[Eq. (37)]{li2014} to construct a test when the number of outliers is known, when each outlier follows the same anomalous distribution and when there is no null hypothesis.

\section{Main Results}
\subsection{Preliminaries}
We first present necessary definitions. Recall the definition of the distribution $P_{\rm{Mix}}^{(\calB,\calC,P_\rmN,\bP_\calB)}$ in \eqref{def:Pmix}. Given any $\calB\in\calS$ and any tuple of distributions $\bP_\calB=(P_\rmN,P_{\rmA,1},\ldots,P_{\rmA,|\calB|})$, define the following information densities (log likelihoods):
\begin{align}
\!\!\imath_{1,l}(x|\calB,\calC,P_\rmN,\bP_\calB)&:=\log\frac{P_{\rmA,l}(x)}{P_{\rm{Mix}}^{(\calB,\calC,P_\rmN,\bP_\calB)}(x)},~l\in[|\calB|]\label{def:i1T},\\
\imath_2(x|\calB,\calC,P_\rmN,\bP_\calB)&:=\log\frac{P_\rmN(x)}{P_{\rm{Mix}}^{(\calB,\calC,P_\rmN,\bP_\calB)}(x)}\label{def:i2T}.
\end{align}
One can verify that $\mathrm{GD}(\calB,\calC,P_\rmN,\bP_\calB)$ in \eqref{def:GDBC} is the linear combination of expectations of the information densities, i.e., 
\begin{align}
\nn&\mathrm{GD}(\calB,\calC,P_\rmN,\bP_\calB)\\*
\nn&:=\sum_{i\in(\calB\cap\calM_{\calC})}\mathbb{E}_{P_{\rmA,\jmath_\calB(i)}}[\imath_{1,\jmath_\calB(i)}(X|\calB,\calC,P_\rmN,\bP_\calB)]\\*
&\qquad+\sum_{i\in(\calM_{\calB}\cap\calM_{\calC})}\mathbb{E}_{P_\rmN}[\imath_2(X|\calB,\calC,P_\rmN,\bP_\calB)].
\end{align}
Furthermore, define the following linear combination of variances of information densities:
\begin{align}
\nn&\rmV(\calB,\calC,P_\rmN,\bP_\calB)\\*
\nn&:=\sum_{i\in(\calB\cap\calM_{\calC})}\mathrm{Var}_{P_{\rmA,\jmath_\calB(i)}}[\imath_{1,\jmath_\calB(i)}(X|\calB,\calC,P_\rmN,\bP_\calB)]\\*
&\qquad+\sum_{i\in(\calM_{\calB}\cap\calM_{\calC})}\mathrm{Var}_{P_\rmN}[\imath_2(X|\calB,\calC,P_\rmN,\bP_\calB)]\label{def:VBC}.
\end{align}
Given any $(\calB,\calC)\in\calS^2$ and any variables $(x_1,\ldots,x_M)$, define the following linear combination of the information densities:
\begin{align}
\nn&\imath_{\calB,\calC}(x_1,\ldots,x_M|P_\rmN,\bP_\calB)
:=\!\!\!\!\sum_{j\in(\calB\cap\calM_\calC)}\!\!\!\!\imath_{1,\jmath_\calB(j)}(x_j|\calB,\calC,P_\rmN,\bP_\calB)\\*
&\qquad+\sum_{\barj\in(\calM_\calB\cap\calM_\calC)}\imath_2(x_{\barj}|\calB,\calC|P_\rmN,\bP_\calB).
\end{align}
Recall that $\calS_\calB$ denotes the set $\calS\setminus\{\calB\}$ and let the elements in $\calS_\calB$ be ordered as $\{\calC_1,\ldots,\calC_{|\calS|-1}\}$. For each $(i,k)\in[|\calS|-1]^2$ such that $i\neq k$, define the following covariance function
\begin{align}
\nn&\mathrm{Cov}(\calC_i,\calC_k,P_\rmN,\bP_\calB):=\mathrm{E}[\imath_{\calB,\calC_i}(X_1,\ldots,X_M|P_\rmN,\bP_\calB)\\*
&\qquad\times\imath_{\calB,\calC_k}(X_1,\ldots,X_M|P_\rmN,\bP_\calB)].
\end{align}
Then define a covariance matrix $\bV(\calB,P_\rmN,\bP_\calB)=\{V_{i,j}(\calB,P_\rmN,\bP_\calB)\}_{(i,j)\in[|\calS|-1]^2}$ where
\begin{align}
\!\!\!V_{i,j}(\calB,P_\rmN,\bP_\calB)
&\!=\left\{
\begin{array}{ll}
\rmV(\calB,\calC_i,P_\rmN,\bP_\calB)&\mathrm{if~}i=j,\\
\mathrm{Cov}(\calC_i,\calC_k,P_\rmN,\bP_\calB)&\mathrm{otherwise.}
\end{array}
\right.
\end{align}
For any $k\in\bbN$, $\rmQ_k(x_1,\ldots,x_k;\bmu,\bSigma)$ is the multivariate generalization of the complementary Gaussian cdf defined as follows:
\begin{align}
\rmQ_k(x_1,\ldots,x_k;\bmu,\bSigma)
&:=\int^{\infty}_{x_1}\ldots\int^{\infty}_{x_k}\calN(\bx;\bmu;\bSigma)\rmd \bx\label{def:kQ},
\end{align}
where $\calN(\bx; \bmu;\bSigma)$ is the pdf of a $k$-variate Gaussian with mean $\bmu$ and covariance matrix $\bSigma$~\cite{tan2015asymptotic}. Furthermore, for any $k\in\bbN$, we use $\mathbf{1}_k$ to denote a row vector of length $k$ with all elements being one and we use $\mathbf{0}_k$ similarly. The complementary cdf $\rmQ_k(\cdot)$ in \eqref{def:kQ}, together with $\mathrm{GD}(\calB,\calC,P_\rmN,\bP_\calB)$ and $\bV(\calB,P_\rmN,\bP_\calB)$, will be critical to upper bound the false reject probabilities.

Finally, given any $\lambda\in\bbR_+$ and any tuple of distributions $\bP_\calB=(P_\rmN,P_{\rmA,1},\ldots,P_{\rmA,T})\in\calP_T(\calX)$, for each $\calB\in\calS$, define the following quantity:
\begin{align}
\mathrm{LD}_{\calB}&(\lambda,P_\rmN,\bP_\calB)
\nn:=\min_{(\calC,\calD)\in\calS^2:\calC\neq\calD}\min_{\substack{\bQ\in(\calP(\calX))^M:\\\rmG_\calC(\bQ)\leq \lambda,~\rmG_\calD(\bQ)\leq\lambda}}\\*
&\Big(\sum_{i\in\calB}D(Q_i\|P_{\rmA,\jmath_\calB(i)})+\sum_{i\in\calM_\calB}D(Q_i\|P_\rmN)\Big)\label{def:LDT}.
\end{align}
The quantity $\mathrm{LD}_{\calB}(\lambda,P_\rmN,\bP_\calB)$ will characterize the false reject exponent under each hypothesis.

\subsection{Second-Order Asymptotics}
We first provide an achievability result, where the performance of the test in \eqref{test:Toutlier} is characterized in terms of misclassification and false alarm probabilities that decay exponentially fast when the false reject probability is upper bounded by a function of the threshold $\lambda$. Furthermore, we demonstrate the optimality of the test in \eqref{test:Toutlier} under the generalized Neyman-Pearson criterion.

\begin{theorem}
\label{result:Toutlier}
For any nominal distribution $P_\rmN$ and anomalous distributions $\bP_T=(P_{\rmA,1},\ldots,P_{\rmA,T})$, given any positive real number $\lambda\in\bbR_+$, the test in \eqref{test:Toutlier} satisfies that for each $\calB\in\calS$,
\begin{align}
\!\!\!\beta_{\calB}(\Psi_n|P_\rmN,\bP_T)&\leq \exp\Big(-n\lambda+|\calX|\log((M-1)n+1)\Big),\\
\!\!\!\rmP_{\mathrm{fa}}(\Psi_n|P_\rmN,\bP_T)\nn&\leq |\calS|^2\exp\Big(-n\lambda\\*
&\qquad\quad+|\calX|\log((M-1)n+1)\Big),
\end{align}
and 
\begin{align}
\nn&\zeta_\calB(\Psi_n|P_\rmN,\bP_T)\leq 1+O\left(\frac{1}{\sqrt{n}}\right)\\*
&\!\!-\rmQ_{|\calS|-1}\big(\sqrt{n}\bar{\mu}(\lambda,P_\rmN,\bP_\calB);\mathbf{0}_{|\calS|-1};\mathbf{0}_{|\calS|-1};\bV(\calB,P_\rmN,\bP_\calB)\big)\label{sr:tout},
\end{align} 
where $\bar{\mu}(\lambda,P_\rmN,\bP_\calB)$ denotes the length-$(|\calS|-1)$ vector with elements $(\lambda-\mathrm{GD}(\calB,\calC_1,P_\rmN,\bP_\calB),\ldots,\lambda-\mathrm{GD}(\calB,\calC_{|\calS|-1},P_\rmN,\bP_\calB))+\mathbf{1}_{|\calS-1|}\times O(\log n/n)$.
\end{theorem}
The proof of Theorem \ref{result:Toutlier} is available in~\cite[Appendix D]{zhou2020second}. Several remarks are as follows.

When the number of outliers $M$ is finite, both misclassification and false alarm probabilities decay exponentially fast, with a speed lower bounded by $\lambda$ asymptotically when $n$ tends to infinity. On the other hand, the false reject under each hypothesis $\rmH_\calB$ is upper bounded by a function of $\lambda$ and critical quantities $\mathrm{GD}(\calB,\calC,P_\rmN,\bP_\calB)$ and $\bV(\calB,P_\rmN,\bP_\calB)$. Note that the threshold $\lambda$ trades off the lower bound on the decay rate of the homogeneous error exponent of the misclassification and false alarm probabilities and the upper bound on the false reject probability. If $\lambda$ increases, the homogeneous error exponent increases while the false reject probability increases as well. Thus better performance in misclassification error and false alarm probabilities leads to worse false reject probabilities.

Asymptotically as $n\to\infty$, if the threshold $\lambda<\min_{i\in[|\calS|-1]}\mathrm{GD}(\calB,\calC_i,P_\rmN,\bP_\calB)=:\mathrm{GD}(\calB,P_\rmN,\bP_{\calB})$, then the false reject probability under hypothesis $\rmH_\calB$ vanishes. One might also be interested in the more practical non-asymptotic case where $n$ is finite. Obtaining the exact solution to such case is almost impossible. However, a second-order asymptotic approximation to the non-asymptotic performance is possible using the result in \eqref{sr:tout}. For this purpose, let $d(\calB)$ be the number of elements in the vector that equals the minimal value, i.e., $d(\calB):=|\{i\in[|\calS|-1]:|\mathrm{GD}(\calB,\calC_i,P_\rmN,\bP_\calB)=\mathrm{GD}(\calB,P_\rmN,\bP_{\calB})\}$, and given any $\varepsilon\in(0,1)$, let
\begin{align}
&L^*(\varepsilon|\calB,P_\rmN,\bP_\calB)
\nn:=\max\Big\{L\in\bbR:\\*
&\quad\rmQ_{d(\calB)}(L\times\mathbf{1}_{d(\calB)};\mathbf{0}_{d(\calB)};\bV(\calB,P_\rmN,\bP_\calB))\geq 1-\varepsilon\Big\}\label{def:L*PT},\\
\nn&\lambda^*(n,\varepsilon|\calB,P_\rmN,\bP_\calB)\\*
&\quad:=\mathrm{GD}(\calB,P_\rmN,\bP_{\calB})+\frac{L^*(\varepsilon|\calB,P_\rmN,\bP_\calB)}{\sqrt{n}}\label{def:lambda^*moutlier}.
\end{align}

We then have the following corollary of Theorem \ref{result:Toutlier}.
\begin{corollary}
\label{second:tout}
For any $(\calB,P_\rmN,\bP_\calB)$, if $\lambda$ satisfies $\lambda\leq \lambda^*(n,\varepsilon|\calB,P_\rmN,\bP_\calB)$ for all $n\in\bbN$, then as the blocklength $n$ increases, the upper bound on the false reject probability under $(\calB,P_\rmN,\bP_\calB)$ tends to $\varepsilon\in(0,1)$, i.e., $\lim_{n\to\infty}\zeta_\calB(\Psi_n|P_\rmN,\bP_T)\leq \varepsilon$.
\end{corollary}

The result in Corollary \ref{second:tout} implies a phase transition phenomenon for our test. In particular, if the threshold $\lambda$ is strictly greater than $\mathrm{GD}_M(P_N,P_A)$, then asymptotically the false reject probabilities tend to one. On the other hand, if $\lambda<\mathrm{GD}_M(P_N,P_A)$, then asymptotically the false reject probabilities vanish. 

Note that $\lambda^*(n,\varepsilon|\calB,P_\rmN,\bP_\calB)$ is a critical bound for the threshold in the test, which trades off a lower bound $\lambda$ on the exponential decay rates of misclassification error and false alarm probabilities and a non-vanishing upper bound $\varepsilon\in(0,1)$ for the maximal false reject probability. Such a result is known as a second-order asymptotic result since it provides a formula for the second dominant term $\frac{L^*(\varepsilon|\calB,P_\rmN,\bP_\calB)}{\sqrt{n}}$ beyond the leading constant term $\mathrm{GD}(\calB,P_\rmN,\bP_{\calB})$ asymptotically as $n\to\infty$.

Finally, we discuss the influence of the number of observed sequences $M$ on the performance of the test in \eqref{test:Toutlier}. As demonstrated in the above remark, $\mathrm{GD}(\calB,P_\rmN,\bP_{\calB})$ is the critical quantity that is related with the performance of the test. Thus, it suffices to study the properties of $\mathrm{GD}(\calB,P_\rmN,\bP_{\calB})$ as a function of $M$ under each hypothesis $\rmH_\calB$. However, it is challenging to obtain closed form equations for the dependence of $\mathrm{GD}(\calB,P_\rmN,\bP_{\calB})$ on $M$ when each outlier is generated from a unique anomalous distributions. Thus, we specialize our results to the case where all anomalous distributions are the same and given by $P_\rmA$. Under this assumption, we have
\begin{align}
\mathrm{GD}(\calB,P_\rmN,\bP_{\calB})
\nn&=\min_{t\in[T]}\min_{l\in[|\calB|]}\big(lD(P_\rmA\|P_{\mathrm{Mix}}^{t,l})\\*
&\qquad+(M-t-l)D(P_\rmN\|P_{\mathrm{Mix}}^{l,t})\big)\label{allsame},
\end{align}
where $P_{\mathrm{mix}}^{t,l}=\frac{lP_\rmA+(M-t-l)P_\rmN}{M-t}$. For any $(t,l)\in[T]\times[|\calB|]$, we have
\begin{align}
\frac{\partial \mathrm{GD}(\calB,P_\rmN,\bP_{\calB})}{\partial M}
&=D(P_\rmN\|P_{\mathrm{Mix}}^{l,t}).
\end{align}
Thus, $\mathrm{GD}(\calB,P_\rmN,\bP_{\calB})$ increases in $M$ if $D(P_\rmN\|P_{\mathrm{Mix}}^{l,t})>0$, which holds almost surely for all distinct nominal and anomalous distributions. This implies that the performance of the test in \eqref{test:Toutlier} increases as the number of observed sequences $M$ increases when the number of outliers $|\calB|$ remains unchanged. On the other hand, the result in \eqref{allsame} implies that for a fixed number of observed sequences $M$, the performance of the test in \eqref{test:Toutlier} degrades as the number of outliers $|\calB|$ increases.

With the above achievability result on the performance of the test in \eqref{test:Toutlier}, it remains to show that the test is in fact optimal in a certain sense. Since nominal and anomalous distributions are unknown, in order to derive a converse result, the classical Neyman-Pearson criterion, which requires knowledge of generating distributions, is not applicable. Furthermore, as proved in \cite{li2014}, for our problem, it is impossible to ensure that all three kinds of error probabilities decay exponentially for all pairs of nominal and anomalous distributions. As a compromise, we adopt the generalized Neyman-Pearson criterion of Gutman~\cite{gutman1989asymptotically} to derive a lower bound on the false reject probability. The generalized Neyman-Pearson criterion is that both misclassification error and false alarm probabilities decay exponentially fast with homogeneous speed for \emph{all} tuple of nominal and anomalous distributions. We give a lower bound on the false reject probability for any particular pair of distributions $(P_\rmN,\bP_T)$ in the following theorem.
\begin{theorem}
\label{converse:m}
Given any $\lambda\in\bbR_+$, for any test $\phi_n$ such that
\begin{align}
\beta_\calB(\phi_n|\tilP_\rmN,\tilde{\bP}_T)\leq \exp(-n\lambda),~\forall~(\tilP_\rmN,\tilde{\bP}_T)\label{Tout:conreqire},
\end{align}
then for any tuple of nominal and anomalous distributions $(P_\rmN,\bP_T)$, for each $\calB\in\calS$, $\zeta_\calB(\Psi_n|P_\rmN,\bP_T)$ is lower bounded by the right hand side of \eqref{sr:tout}.
\end{theorem}
The proof of Theorem \ref{converse:m} is available in \cite[Appendix E]{zhou2020second}. The result in Theorem \ref{converse:m} holds for any number of observed sequences $M$ and when the length $n$ of each observed sequence $n$ is such that $O(\frac{\log n}{n})$ and $O(\frac{1}{\sqrt{n}})$ can be neglected. 

\subsection{Large Deviations}
We next characterize the tradeoff between the false reject exponent and the homogeneous error exponent of the misclassification error and false alarm probabilities. Recall the definition of $\mathrm{LD}_{\calB}(\lambda,P_\rmN,\bP_\calB)$ in \eqref{def:LDT}.
\begin{theorem}
\label{result:Tout:exp}
For any nominal distribution $P_\rmN$ and anomalous distributions $\bP_T=(P_{\rmA,1},\ldots,P_{\rmA,T})$, given any positive real number $\lambda\in\bbR_+$, the test in \eqref{test:Toutlier} satisfies that for each $\calB\in\calS$,
\begin{align}
\liminf_{n\to\infty}-\frac{1}{n}\log\beta_{\calB}(\Psi_n|P_\rmN,\bP_T)&\geq\lambda,\\
\liminf_{n\to\infty}-\frac{1}{n}\log\rmP_{\mathrm{fa}}(\Psi_n|P_\rmN,\bP_T)&\geq \lambda,\\
\liminf_{n\to\infty}-\frac{1}{n}\log\zeta_\calB(\Psi_n|P_\rmN,\bP_T)&\geq\mathrm{LD}_{\calB}(\lambda,P_\rmN,\bP_\calB).
\end{align} 
Conversely, for any test that ensures the homogeneous exponential decay rate of the misclassification error and false alarm is no less than $\lambda$ for all tuples of nominal and anomalous distributions, under any nominal distribution $P_\rmN$ and anomalous distributions $\bP_T=(P_{\rmA,1},\ldots,P_{\rmA,T})$, the false reject exponent is also upper bounded by $\mathrm{LD}_{\calB}(\lambda,P_\rmN,\bP_\calB)$ under each hypothesis $\rmH_\calB$.
\end{theorem}
The proof idea of Theorem \ref{result:Tout:exp} is available in \cite[Section III.E]{zhou2020second}. We note that the threshold $\lambda$ governs the tradeoff between the false reject exponent and the homogeneous error exponent under each hypothesis. From the definition of $\mathrm{LD}_{\calB}(\lambda,P_\rmN,\bP_\calB)$ in \eqref{def:LDT}, it follows that the false reject exponent $\mathrm{LD}_{\calB}(\lambda,P_\rmN,\bP_\calB)$ decreases in $\lambda$. One can show that $\mathrm{LD}_{\calB}(\lambda,P_\rmN,\bP_\calB)>0$ if and only if $\lambda<\mathrm{GD}(\calB,P_\rmN,\bP_\calB)$ and the maximal false reject exponent satisfies
\begin{align}
\nn&\max_{\lambda\in(0,\mathrm{GD}(\calB,P_\rmN,\bP_\calB))}\mathrm{LD}_{\calB}(\lambda,P_\rmN,\bP_\calB)\\*
&\leq \min_{Q\in\calP(\calX)} \Big(\sum_{i\in\calB}D(Q\|P_{\rmA,\jmath_{\calB}(i)})+(M-|\calB|)D(Q\|P_\rmN)\Big).
\end{align}
Therefore, if the threshold $\lambda<\min_{\calB\in\calS}\mathrm{GD}(\calB,P_\rmN,\bP_\calB)$, then regardless of the number of outliers, the misclassification error, the false alarm and false reject probabilities decay exponentially fast for any tuple of distributions $(P_\rmN,\bP_T)$ such that $\min_{\calB\in\calS}\mathrm{GD}(\calB,P_\rmN,\bP_\calB)$ is strictly positive.

Note that asymptotically, the exponents of probabilities of misclassification error and false alarm are equal. This is an artifact of our test in \eqref{test:Toutlier} where only one threshold $\lambda$ is used. It would be worthwhile to investigate tests that can fully characterize the exponent tradeoff of all three kinds of error probabilities, beyond the degenerate ``corner-point'' case in this paper. Such investigations will be pursued in future work.

\section{Conclusion}
\label{sec:conclude}
We revisited the outlier hypothesis testing problem studied by Li \emph{et al.}~\cite{li2014} and proposed a threshold-based test, analyzed its performance and proved its asymptotic optimality under the generalized Neyman-Pearson criterion~\cite{gutman1989asymptotically}. In the future, one can consider a sequential setting~\cite{li2017universal} and study sequences with continuous~\cite{Zou2017TSP} or large alphabets~\cite{kelly2013}.

\newpage
\bibliographystyle{IEEEtran}
\bibliography{IEEEfull_lin}
\end{document}

%% file: preamble.tex
\usepackage[mathscr]{eucal}
\usepackage[cmex10]{amsmath}
\usepackage{epsfig,epsf,psfrag}
\usepackage{amssymb,amsmath,amsthm,amsfonts,latexsym}
\usepackage{amsmath,graphicx,bm,xcolor,url,overpic}
\usepackage{fixltx2e}
\usepackage{array}
\usepackage{verbatim}
\usepackage{bm}
\usepackage{algorithmic}
\usepackage{algorithm}
\usepackage{verbatim}
\usepackage{textcomp}
\usepackage{mathrsfs}
\usepackage{epstopdf}

\newcommand{\openone}{\leavevmode\hbox{\small1\normalsize\kern-.33em1}}

\catcode`~=11 \def\UrlSpecials{\do\~{\kern -.15em\lower .7ex\hbox{~}\kern .04em}} \catcode`~=13 

\allowdisplaybreaks[4]

\newcommand{\nn}{\nonumber}

\newcommand{\calB}{\mathcal{B}}
\newcommand{\calC}{\mathcal{C}}
\newcommand{\calD}{\mathcal{D}}

\newcommand{\calM}{\mathcal{M}}
\newcommand{\calN}{\mathcal{N}}

\newcommand{\calP}{\mathcal{P}}

\newcommand{\calS}{\mathcal{S}}

\newcommand{\calX}{\mathcal{X}}

\newcommand{\bP}{\mathbf{P}}

\newcommand{\bQ}{\mathbf{Q}}

\newcommand{\bV}{\mathbf{V}}

\newcommand{\bx}{\mathbf{x}}
\newcommand{\bX}{\mathbf{X}}

\newcommand{\rmA}{\mathrm{A}}

\newcommand{\rmd}{\mathrm{d}}

\newcommand{\rmG}{\mathrm{G}}

\newcommand{\rmH}{\mathrm{H}}

\newcommand{\rmN}{\mathrm{N}}

\newcommand{\rmP}{\mathrm{P}}

\newcommand{\rmQ}{\mathrm{Q}}
\newcommand{\rmr}{\mathrm{r}}

\newcommand{\rmS}{\mathrm{S}}

\newcommand{\rmT}{\mathrm{T}}

\newcommand{\rmV}{\mathrm{V}}

\newcommand{\bbN}{\mathbb{N}}

\newcommand{\bbP}{\mathbb{P}}

\newcommand{\bbR}{\mathbb{R}}

\DeclareMathAlphabet{\mathbsf}{OT1}{cmss}{bx}{n}
\DeclareMathAlphabet{\mathssf}{OT1}{cmss}{m}{sl}

\DeclareSymbolFont{bsfletters}{OT1}{cmss}{bx}{n}  
\DeclareSymbolFont{ssfletters}{OT1}{cmss}{m}{n}
\DeclareMathSymbol{\bsfGamma}{0}{bsfletters}{'000}
\DeclareMathSymbol{\ssfGamma}{0}{ssfletters}{'000}
\DeclareMathSymbol{\bsfDelta}{0}{bsfletters}{'001}
\DeclareMathSymbol{\ssfDelta}{0}{ssfletters}{'001}
\DeclareMathSymbol{\bsfTheta}{0}{bsfletters}{'002}
\DeclareMathSymbol{\ssfTheta}{0}{ssfletters}{'002}
\DeclareMathSymbol{\bsfLambda}{0}{bsfletters}{'003}
\DeclareMathSymbol{\ssfLambda}{0}{ssfletters}{'003}
\DeclareMathSymbol{\bsfXi}{0}{bsfletters}{'004}
\DeclareMathSymbol{\ssfXi}{0}{ssfletters}{'004}
\DeclareMathSymbol{\bsfPi}{0}{bsfletters}{'005}
\DeclareMathSymbol{\ssfPi}{0}{ssfletters}{'005}
\DeclareMathSymbol{\bsfSigma}{0}{bsfletters}{'006}
\DeclareMathSymbol{\ssfSigma}{0}{ssfletters}{'006}
\DeclareMathSymbol{\bsfUpsilon}{0}{bsfletters}{'007}
\DeclareMathSymbol{\ssfUpsilon}{0}{ssfletters}{'007}
\DeclareMathSymbol{\bsfPhi}{0}{bsfletters}{'010}
\DeclareMathSymbol{\ssfPhi}{0}{ssfletters}{'010}
\DeclareMathSymbol{\bsfPsi}{0}{bsfletters}{'011}
\DeclareMathSymbol{\ssfPsi}{0}{ssfletters}{'011}
\DeclareMathSymbol{\bsfOmega}{0}{bsfletters}{'012}
\DeclareMathSymbol{\ssfOmega}{0}{ssfletters}{'012}

\newcommand{\tilP}{\tilde{P}}

\newcommand{\hatT}{\hat{T}}

\newcommand{\barj}{\bar{j}}

\newcommand{\bmu}{\bm{\mu}}

\newcommand{\bSigma	}{\bm{\Sigma}}

\newtheorem{theorem}{Theorem}

\newtheorem{corollary}{Corollary}